\documentstyle[12pt]{article}

\def\today{9.10.03}

\newcommand\equ[1]{{\rm (\ref{#1})}}
\newcommand\beq[1]{ \begin{equation}\label{#1} }
\newcommand{\eeq}{ \end{equation} }

\newcommand\beqa[1]{ \begin{eqnarray} \label{#1}}
\newcommand{\eeqa}{ \end{eqnarray} }
\newcommand{\beqano}{ \begin{eqnarray*} }
\newcommand{\eeqano}{ \end{eqnarray*} }

\newtheorem{theorem}{Theorem}[section]
\newtheorem{definition}{Definition}[section]
\newtheorem{proposition}{Proposition}[section]
\newtheorem{lemma}{Lemma}[section]
\newtheorem{remark}{Remark}[section]

\newtheorem{corollary}{Corollary}[section]

\newcommand\dfn[1]{ \begin{definition}\label{#1} }
\newcommand\thm[1]{ \begin{theorem}\label{#1}}

\newcommand\thmtwo[2]{ \begin{theorem}[#2]\label{#1}}
\newcommand\ethm{ \end{theorem} }
\newcommand\pro[1]{ \begin{proposition}\label{#1}}
\newcommand\protwo[2]{ \begin{proposition}[#2]\label{#1}}
\newcommand\epro{ \end{proposition} }
\newcommand\lem[1]{ \begin{lemma}\label{#1}}
\newcommand\lemtwo[2]{ \begin{lemma}[#2]\label{#1}}
\newcommand\elem{ \end{lemma} }

\newcommand{\giu}{{\medskip\noindent}}

\newcommand{\torus}{ {\bf T}   }\renewcommand{\natural}{ {\bf N}   }
\newcommand{\real}{ {\bf R}   }
\newcommand{\integer}{ {\bf Z}   }
\newcommand{\complex}{ {\bf C}   }

\renewcommand{\d}{ {\delta}   }

\newcommand{\z}{ {\zeta} }

\newcommand{\cF }{ {\cal F} }
\newcommand{\cH}{ {\cal H} }
\newcommand{\cU}{ {\cal U} }
\newcommand{\cI}{ {\cal I} }

\newcommand{\cL}{ {\cal L} }

\newcommand{\cP}{ {\cal P} }

\newcommand{\cS}{ {\cal S} }
\newcommand{\cT}{ {\cal T} }

\newcommand{\ov}{\overline}

\newcommand{\bz}{ {\overline z} }

\renewcommand{\o}{ {\omega}   }
\renewcommand{\O}{ {\Omega}   }

\newcommand{\ii}{ {\rm i} }

\begin{document}

\author{Dario Bambusi\footnote{
Dipartimento di Matematica, Via Saldini 50, 20133, Milano, Italy.   
Supported by GNFM.}, Massimiliano Berti\footnote{
SISSA, Via Beirut 2-4, 34014, Trieste, Italy.
Supported by MURST, {\sl Variational Methods and Nonlinear
Differential Equations}.}}



\title{A BIRKHOFF--LEWIS TYPE THEOREM FOR SOME
  HAMILTONIAN PDEs}

\date{\today}


\maketitle

\begin{abstract}
In this paper we give an extension of the Birkhoff--Lewis theorem to
some semilinear PDEs. Accordingly we prove existence of infinitely
many periodic orbits with large period accumulating at the
origin. Such periodic orbits bifurcate from resonant finite
dimensional invariant tori of the fourth order normal form of the
system. Besides standard nonresonance and nondegeneracy assumptions,
our main result is obtained assuming a regularizing property of the
nonlinearity. We apply our main theorem to a semilinear beam equation
and to a nonlinear Schr\"odinger equation with smoothing nonlinearity.
\end{abstract}

\setcounter{section}{0}
\section{Introduction}
\label{s1}

In 1934 Birkhoff and Lewis [BL34] (see also [Lew34, Mos77]) proved their
celebrated theorem on existence of periodic orbits with large period
close to elliptic equilibria of Hamiltonian
systems\footnote{Actually [BL34] considers the
neighborhood of an elliptic, non-constant, periodic orbit, but the scheme is
essentially the same for elliptic equilibria.}. Here we give a
generalization of their result to some semilinear Hamiltonian
PDEs. 

Birkhoff--Lewis procedure consists in putting the system in
fourth order (Birkhoff) normal form, namely in the form
\begin{equation}
\label{z4}
H=H_0+G_4+R_5, \, \quad H_0:=\sum_{j= 1}^n\omega_j\frac{p_j^2+q_j^2}{2},
\end{equation}
where $G_4$ is a homogeneous polynomial of degree 2 in the actions
$I_j:=(p_j^2+q_j^2)/2$ and $R_5$ is a remainder having a zero of fifth
order at the origin. Then system (\ref{z4}) is a perturbation of the
integrable system $H_0+G_4$. Under a nondegeneracy condition (that
also plays a fundamental role in KAM theory) the action to frequency
map of this integrable system is one to one, and therefore there exist
infinitely many resonant tori on which the motion is periodic. The
question is: Do some of these periodic orbits persist under the
perturbation due to the term $ R_5 $?  Birkhoff--Lewis used the
implicit function theorem and a
topological argument to prove that there exists a sequence of
resonant tori accumulating at the origin with the property that at
least two periodic orbits bifurcate from each one of them.

In order to extend this result to infinite dimensional systems
describing Hamiltonian PDEs one meets two difficulties: the first one
is the generalization of Birkhoff normal form to PDEs and the second
one is the appearing of a small denominator problem.

Here we decide to work in a way which is as straightforward as
possible, so, instead of considering the standard Birkhoff normal form
of the system, whose extension to PDEs is not completely understood at
present\footnote{See however the recent works [Bam03,BG03].}, we
consider its ``seminormal form'', namely the kind of normal form
employed to construct lower dimensional tori.  Precisely, having fixed
a positive $n$, we split the phase variables in two groups, namely the
variables with index smaller than $n$ and the variables with index
larger than $n$. We will denote by $\hat z$ the whole set of variables
with index larger than $n$. We construct a canonical transformation
putting the system in the form
\begin{equation}
\label{semi}
H_0+\overline G+\hat G+K
\end{equation}
where $\overline G$ depends only on the actions, $\hat G$ is at least
cubic in the variables $\hat z$ with index larger than $n$, and $K$
has a zero of sixth order at the origin.  The interest of such a
seminormal form is that the normalized system $H_0+\overline G+\hat G$
has the invariant $2n$--dimensional manifold $\hat z=0$ which is
filled by $n$--dimensional invariant tori.  

Under a nondegeneracy condition, the frequencies of the flow in such
tori cover an open subset of $ \real^n $.  We concentrate on the
resonant tori filled by periodic orbits and we prove that at least $n$
geometrically distinct periodic orbits of each torus survive the
perturbation due to the term $K$. Since the orbits bifurcate from
lower dimensional tori we have to impose a further nondegeneracy
condition in order to avoid resonances between the frequency of the
periodic orbit and the frequencies of the transversal oscillations.

The proof is based on a variational Lyapunov Schmidt reduction similar
to that employed in [BBV03] and inspired to [ACE87]. It turns out that in
the present case the range equation involves small denominators. To
solve the corresponding problem we use an approach similar to that of
[Bam00]. In particular we impose a strong condition on the small
denominators and we show that, if the vector field of the nonlinearity
is smoothing, then the range equation can be solved by the contraction
mapping principle.  Next, the kernel equation is solved by noting that
it is the Euler-Lagrange equation of the action functional restricted
to the solutions of the range equation. The restricted functional
turns out to be defined on $ \torus^n $ and so existence and
multiplicity of solutions (critical points) follows by the classical
Lusternik-Schnirelmann theory.

Finally we apply the general theorem to 
the nonlinear beam equation 
\begin{equation}
\label{beam}
u_{tt}+u_{xxxx}+mu=f(u)\ .
\end{equation}
with Dirichlet boundary conditions on a segment. We consider $m$ as a
parameter varying in the segment $[0,L]$ and we show that the
assumptions of the abstract theorem are fulfilled provided one
excludes from the interval a finite number of values of $m$.  As a
second application we will deal with a nonlinear Schr\"o\-dinger
equation with a smoothing nonlinearity of the type considered in
[P\"os02].

We recall that families of periodic solutions to Hamiltonian PDEs have
been constructed by many authors (see e.g. [Kuk93m, Way90, P\"os96b,
CW93, Bou95, Bam00]). The main difference is that the periodic orbits
of the above quoted papers are a continuation of the linear normal
modes to the nonlinear system. In particular their period is close to
one of the periods of the linearized system. Moreover (except in the
resonant case, see [BP01,BB03]) each periodic solution involve only
one of the linear oscillators\footnote{in the sense that all the other
ones have a much smaller amplitude of oscillation.}.

On the contrary, the periodic orbits constructed in the present paper
are the shadows of resonant tori; they are a purely nonlinear
phenomenon, have long period, and moreover each periodic motion
involves $n$ linear oscillators that oscillate with amplitudes of the
same order of magnitude. 

\bigskip

\noindent
{{\bf Acknowledgements:} {\it The authors thanks L. Biasco,
P. Bolle and M. Procesi for interesting discussions.}}

\section{Main result}
\label{main1}

Consider a real Hamiltonian system with real\footnote{i.e.
if $\overline z$ is actually the complex conjugate to $z$ then the
  Hamiltonian $H$ takes real values.} Hamiltonian function
\begin{equation}
\label{hami}
H(z,\overline z)=\sum_{j\geq 1}\omega_jz_j\overline 
z_j+P(z,\overline z) \equiv H_0+P
\end{equation}
where $P$ has a zero of third order at the origin and the symplectic
structure is given by $\ii \sum_j dz_j \wedge d\overline z_j$. Here
$z$ and $\bar z$ are considered as independent variables. Often we
will write only the equation for $z$ since the equation for $\bar z$
is obtained by complex conjugation.
The formal 
Hamiltonian vector field of the system is $X_{H}(z,\overline z):=(
\ii \frac{\partial H}{\partial \overline z_j}, - \ii\frac{\partial
H}{\partial z_j})$, and
therefore the equations of motion have the form
\begin{equation}
\label{equa}
\dot z_j= \ii \omega_j z_j+ \ii 
\frac{\partial P}{\partial \overline z_j}\ ,\quad
\dot {\overline z}_j= 
- \ii \omega_j\overline z_j- \ii \frac{\partial P}{\partial z_j}. 
\end{equation}
Define the complex Hilbert space
$$
{\cal H}^{a,s}( \complex )  
:= \Big\{ w = (w_1, w_2, \ldots ) \in \complex ^\infty \ \Big| \
||w||_{a,s}^2 := 
\sum_{j \geq 1}|w_j|^2 j^{2s} e^{2 j a}  < \infty \Big\}.
$$ 
 We fix $s \geq 0 $ and $a \geq 0 $ and we will 
study the system in the  phase space
$$
\cP_{a,s}:={\cal H}^{a,s}( \complex )  \times {\cal H}^{a,s}(
\complex )\ni (z,\bar z)\ .   
$$

Fix any finite integer $ n \geq 2 $ and denote
$\o := (\o_1, \ldots , \o_n )$, $ \Omega := (\o_{n+1}, \o_{n+2},
\ldots ) $. We assume
\begin{itemize}
\item[(A)] The frequencies grow at least linearly at infinity, namely
  there exists $a>0$ and $d_1\geq 1$ such that 
$$
\omega_j\sim aj^{d_1}.
$$
\item[(NR)] For any
$k\in\integer^n$, $l\in\integer^\infty$ with $|l|\leq 2$ and $|k|+|l|\leq 5$,
one has
\begin{equation}
\label{nonris}
\omega\cdot k+\Omega\cdot l\not=0.
\end{equation} 
\item[(S)] There exists a neighbourhood of the origin $\cU\subset
  \cP_{a,s}$ and $d\geq 0$ such that $X_{P}\in
  C^{\infty}(\cU,\cP_{a,s+d})$.
\end{itemize}

\begin{remark}
\label{ss}
In applications to PDEs, property (S) is usually a consequence of the
smoothness of the Nemitsky operator defined by the nonlinear part of
the equation. In order to ensure (S) one has usually to restrict to
the case where the functions with Fourier coefficients in $\cP_{a,s}$
is an algebra (with the product of convolution between sequences). 
This imposes some limitations to the choice of the
indices $a,s$.
\end{remark}

\begin{proposition} 
\label{Birknf} 
Assume (A,NR,S). There exists a real analytic, symplectic change of
variables $ \cT $ defined in some neighborhood $\cU'\subset\cP_{a,s}$
of the origin, transforming the Hamiltonian $ H $ in seminormal form
up to order six, namely into 
\beq{HBir} 
H\circ \cT\equiv { \cal H } =
H_0 + \overline G + \hat G + K \eeq with
$$
{\ov G} = \frac{1}{2} \sum_{\min (i,j) \leq n} {\ov G}_{ij}  
| z_i |^2 | z_j |^2,  
$$ 
$ {\ov G}_{ij} = {\ov G}_{ji} $,
$ \hat G = O ( || \hat z ||_{a, s}^3 ) $ 
where $ \hat z := ( z_{n+1} , z_{n+2}, 
\ldots )$ and $ K = O(|| z ||_{a, s }^6) $. 
Moreover 
\beq{sreg}
X_{\ov G}, X_{ \hat G }, X_K\in C^{\infty }(\cU',\cP_{a,s+d})\ ,\quad
\|z-\cT(z)\|_{a,s+d}\leq C \|z\|_{a,s}^2\ .
\eeq
\end{proposition}
We defer the proof of this proposition to the Appendix.

The interest of such a seminormal form is that the system
obtained by neglecting the reminder $K$ has the invariant manifold $\hat
z=0$ on which the system is integrable.

As a variant with respect to the standard finite dimensional procedure
we have left the third order term $ \hat G $ but  
normalized the system up to order six (instead of five).
This is needed in lemma \ref{contraction}.

We rewrite the Hamiltonian ${ \cal H }$  in the form  
\beq{Hcal}
{\cal H} := \o \cdot I + \Omega \cdot Z + 
\frac{1}{2} A I \cdot I +
(BI, Z) + \hat G  + K
\eeq
where 
$ I := $ $(|z_1|^2, \ldots, |z_n|^2 ) $, 
$ Z  :=$ $(|z_{n+1}|^2, |z_{n+2}|^2 , \ldots )$ are the actions, 
$ A $ is the $ n \times n  $ matrix  
\begin{equation}
\label{defA}
A =  ({\ov G}_{ij})_{1 \leq i,j \leq n} 
\end{equation}
and $ B $ is the $ \infty \times n $ matrix  
\begin{equation}
B = ({\ov G}_{ij})_{1 \leq j \leq n < i}  
\end{equation}
\begin{remark}
\label{smo1}
Due to (\ref{sreg},\ref{Hcal}), one has  
$|(BI)_j|\leq C|I|j^{-d}$ for a suitable
$C$. Indeed, since   
$ X_{\ov G} $ maps $\cP_{a,s}$ to $\cP_{a,s+d} $,  
the operator $ z_j \mapsto
(BI)_jz_j $ maps $\cP_{a,s}$ to $\cP_{a,s+d}$ and therefore its
eigenvalues  $(BI)_j$ must fulfill the above property.
\end{remark}

Introduce action angle variables for the first $n$ modes  
by $ z_j = | z_j | e^{ \ii \phi_j } =$ $ \sqrt{I_j} e^{ \ii \phi_j } 
$ for $j =1, \ldots, n $.

Perform the rescaling $ I_j \to \eta^2 I_j $, $ \phi_j \to \phi_j $
for $ j = 1, \ldots, n$, $ z_j \to \eta z_j $, 
$ \overline{z}_j \to \eta \overline{z}_j $ for $ j \geq n+1 $ and
divide the Hamiltonian by $ \eta^2 $. We get 
\beq{ham (3)} {\cal H} (
I, \phi , \hat z , \hat{\overline z}) = \o \cdot I+ \Omega \cdot Z +
\eta \hat G_\eta + 
\eta^2 \Big( \frac{1}{2} A I \cdot I + (B I, Z) \Big) + \eta^4 K_\eta.  
\eeq 
where $ \hat G_\eta = O( || \hat z ||_{a,s}^3 )$ and 
$ K_\eta (z) = O(|| z ||_{a,s}^6) $. We will still denote
by $ \cP_{a,s}\equiv \real^n \times$ $\torus^n \times$ ${\cal H}^{a,s}
\times$ ${\cal H}^{a,s}$ the phase space.

We will find periodic solutions of the Hamiltonian system
\equ{ham (3)} close to periodic solutions of the integrable Hamiltonian
system
\beq{HSint}
\dot I = 0, \quad
\dot \phi = \o + \eta^2 ( AI + B^T Z ), \quad
\dot z_j =  \ii \Big( \O_j + \eta^2 ( B I )_j \Big) z_j, \ \ j \geq n+1  
\eeq
in which $\hat G_\eta$ and $K_\eta$ are neglected.
The manifold $ \{ \hat z = 0 \} $ is invariant for the Hamiltonian system
\equ{HSint} and it is completely filled up by the invariant tori
$$
{\cal T}(I_0) := \{ I = I_0, \ \phi \in \torus^n , \ \hat z = 0 \} 
$$
on which the motion is linear with frequencies 
$$
\widetilde \omega\equiv \widetilde\omega(I_0):=\o + \eta^2 A I_0.
$$
Such a torus is linearly stable and the
frequencies of small oscillation about the torus $\cT(I_0)$ are the  
``shifted elliptic frequencies'', namely
\begin{equation}\label{omega eta}
\widetilde \O_{j} (I_0) := (\O +\eta^2 B I_0)_j.
\end{equation}

If all the $\widetilde\omega$'s are integer multiples of a single
frequency, namely if 
\beq{otilde} 
\widetilde \o := \o+\eta^2 A I_0 =
\frac{1}{T} 2 \pi k \in \frac{1}{T} 2 \pi \integer^n, 
\eeq 
then $
{\cal T}(I_0) $ is a {\it completely resonant} torus, supporting the
family of $T$-periodic motions 
\beq{unpert} 
{\cal P} := \Big\{ I(t) =
I_0, \quad \phi (t) = \phi_0 + \widetilde \o t, \quad \hat z(t) = 0
\Big\}.  
\eeq 
The whole family $ { \cal P} $ will not persist in the
dynamics of the complete Hamiltonian system \equ{ham (3)}. We will show
that, under suitable assumptions, at least $ n $
geometrically distinct $T$-periodic solutions persist.
More precisely, we will show that this happens for $\eta$ small enough
and for any choice of $I_0$
and $T$ with 
\begin{equation}
\label{periodi}
\|I_0\|\leq C\ ,\qquad\frac{1}{\eta^2}\leq T\leq \frac{2}{\eta^2}
\end{equation}
where $C$ is independent of $\eta$,
fulfilling
\begin{itemize}
\item[(H1)] Equation (\ref{otilde}) holds. 
\item[(H2)] There exist $\delta>0$ and $\tau<d$ such that 
\beq{Mmeno1}
| \widetilde \O_{j} T - 2 \pi l | \geq \frac{\delta}{j^\tau}, \ \
 \forall l \in \integer,
\ \forall j \geq n + 1. 
\eeq 
\end{itemize}
 
\begin{proposition}
\label{nonresonance}
Fix $\tau>1$. 
Assume (A), ${\rm \det} A\not=0$ and 
\begin{equation}
\label{defhat}
\hat\Omega_j:=\left(\Omega-BA^{-1}\omega\right)_j\not=0 \ ,\quad
\forall j \geq n+1,
\end{equation}
then, for any $\eta>0$, and almost any $T $ fulfilling
(\ref{periodi}) there exists $I_0$ such that (H1,H2) hold.
\end{proposition}
\noindent {\bf Proof.} 
Fix $\eta$, we define $ I_0 := I_0 (T)$ as a function of $T$ so that that
\equ{otilde} is identically satisfied. Then we find $ T $ so that
the non resonance property \equ{Mmeno1} holds. Fix $\eta$
and define, 
\beqa{I0}
I_0 &:=& I_0(T)
:= \frac{2\pi}{\eta^2 T} A^{-1} 
\Big( \Big[ \frac{\o T}{2 \pi} \Big]  - \frac{\o T}{2 \pi} \Big)
,\\ k &:=& k(T) = \Big[ \frac{\o T}{2\pi} \Big] ,\label{kT} 
\eeqa
where $[(x_1,\dots,x_n)] := ( [x_1],\dots, [ x_n ])$ and $[x ] \in
\integer $ denotes the integer part of $x \in \real $.  With the
choice \equ{I0},\equ{kT}, $ \o T + T \eta^2 A I_0 = 2 \pi k $, and
$I_0$ is of order 1 since $ T \eta^2 \geq 1 $. 

We come to the nonresonance property (\ref{Mmeno1}). To study it
remark that the function $ T \to [\o T \slash 2 \pi ] $ is piecewise
constant.  Hence, for any $T_0 \in (\eta^{-2}, 2 \eta^{-2}) $ there
exists an interval $\cI_0= (T_0-a, T_0 + b) \subset [\eta^{-2}, 2
\eta^{-2} ]$ such that $ [\o T \slash 2 \pi ] := k_0 $ is constant for
$T\in\cI_0$.  Moreover, the union of such intervals cover the whole
set of values in which we are interested.  We will construct a subset
of full measure of $\cI_0$, in which condition (H2) is fulfilled.

So, for fixed $j,l$ consider the set
\begin{equation}
\label{resojl}
B_{jl}(\tau,\delta):=\left\{T\in\cI_0\ :\ |\widetilde \Omega_j T-2\pi
l|<\delta/j^{\tau}  \right\}.
\end{equation}
Remark that 
$$
\widetilde \Omega_j T= \hat \Omega_j T+\left(2\pi BA^{-1}
\left[\frac{\omega T}{2\pi} \right] \right)_j,
$$
so that, in $\cI_0$ 
$$
\frac{d}{dT}\left(\widetilde \Omega_j T-2\pi
l\right)=\hat \Omega_j. 
$$
By (A) and remark \ref{smo1} there exists $C$ such that
\begin{equation}
\label{growth2}
\left|\widetilde \Omega_j\right|\geq Cj^{d_1}\ ,\quad \left|\hat
\Omega_j\right|\geq Cj^{d_1}.
\end{equation}
Then $B_{jl}$ is an interval with
length $|B_{jl}|$ controlled by
\begin{equation}
\label{length}
\left|B_{jl}\right|< 2\frac{\delta}{C j^{\tau+d_1}}.
\end{equation}
Fix $j$ and estimate the number of $l$ for which the set $B_{jl}$ is
(possibly) non empty. First remark that, due to (A), one
has that, as $T$ varies in $\cI_0$, the quantity $\widetilde \Omega_jT$
varies in a segment of length smaller than $Cj^{d_1}$, with a suitable
$C$. This means that there are at most $C j^{d_1}$ values of $l$ which
fall in such an interval (with redefined $C$). So, one has
\begin{equation}
\label{sl}
\left|\bigcup_{l}B_{jl}\right|\leq \frac{C\delta}{j^\tau}\ .
\end{equation}
Thus, provided $\tau>1$ as we assumed, one has that
\begin{equation}
\label{mea1}
\left|\bigcup_{jl}B_{jl}\right|\leq C\delta.
\end{equation}
By this estimate, the intersection over $\delta$ of such sets has zero
measure, which is the thesis. $\hfill\Box$

\begin{theorem}
\label{main?}
Let $T$ and $I_0$ fulfill (\ref{periodi}) and (H1,H2), 
then, provided $ \eta $ is small
enough, there exist $n$ geometrically distinct periodic orbits of 
the Hamiltonian system  ${\cal H} $ cfr. (\ref{ham (3)}) with
period $T$ which are $\eta^{2}$ close in $\cP_{a,s}$ to the torus
$\cT(I_0)$.
\end{theorem}

Going back to the original system one has

\begin{corollary}
\label{main}
Consider the Hamiltonian system (\ref{equa}) and fix a positive
$n$. Assume that (A,NR,S) hold, that $ {\rm det} A\not=0$ (cf. (\ref{defA}))
and that $\hat \Omega_j\not=0$ for all $j\geq n+1$
(cf. (\ref{defhat})). Finally assume $d>1$.  \\ Then, for any positive
$\eta\ll1$ there exist at least $n$ distinct periodic orbits
$z^{(1)}(t),...,z^{(n)}(t)$ with the following properties:
\begin{itemize}
\item $\|z^{(l)}(t)\|_{a,s} \leq C\eta$ for $l=1,...,n$ and
  $t\in\real$;

\item $\|\Pi_{>n}z^{(l)}(t)\|_{a,s}\leq C\eta^2$ for $l=1,...,n$ and
  $t\in\real$; here $\Pi_{>n}$ is the projector on the modes with
  index larger than $n$;

\item The period $T$ of $z^{(l)}$ does not depend on $l$ and fulfills
$\displaystyle{ \eta^{-2}\leq T\leq 2\eta^{-2}}$.
\end{itemize}
\end{corollary}

\begin{remark}
\label{fund}
If the integer numbers $( k_1, ... , k_n) = \widetilde \omega T / 2 \pi $ are
relatively prime then $T$ is the minimal period
of the periodic solutions $z^{(l)}$. 
Indeed $z^{(l)}$ are $T$-periodic functions close to the functions defined in 
(\ref{unpert}) which have minimal period $T$.
\end{remark}

\section{Proof of theorem \ref{main?}}

Since the problem is Hamiltonian, any periodic solution of
the system is a critical point of the action functional
\begin{equation}
\label{functional}
S(I,\phi, \hat z , \hat {\overline z})=\int_0^{{T}}\left(I\cdot \dot \phi+
\ii \sum_{j \geq n+1}z_j\dot {\overline z}_j-\cH (I,\phi, \hat z, 
\hat{\overline z})\right)dt
\end{equation}
in the space of $T$-periodic, $\cP_{a,s}$--valued functions. Here
$\cH$ is given by (\ref{ham (3)}). 

We look for a periodic solution $ \z := (\phi, I,\hat z, 
\hat {\overline z} ) $ of the
form 
\beq{variation} \phi (t) = \phi_0 +\widetilde\o t + \psi(t),
\quad I ( t ) = I_0 + J ( t ), \quad 
\eeq 
where $ (\psi, J, \hat z, \hat {\overline z} )$ are periodic functions of 
period $ T $ taking
values in the covering space $ {\bf R}^n \times {\bf R}^n \times {\cal
H}^{a,s} \times {\cal H}^{a,s} $ of $ { \cal P }_{a,s} $ (that for
simplicity will still be denoted by $ { \cal P }_{a,s} $).  Hence ($
\psi $, $J$, $ \hat z $) must satisfy (in the sequel for simplicity 
of notation we will only consider the equation for $\hat z$)
\begin{eqnarray}
\label{sistri}
\nabla_\phi S(\zeta)=0 &\iff & \qquad\qquad \dot J = R_\phi ( \z )
\label{ult} 
\\
\label{31}
\nabla_I S(\zeta)=0 &\iff & \qquad \dot \psi - \eta^2 A \,J  =  R_I ( \z )
\\
\label{32}
\nabla_{\overline z_j} S(\zeta)=0 &\iff & \qquad \dot z_j - \ii \widetilde \O_j
z_j = (R_{\overline z})_j ( \z )
\end{eqnarray}
where
\begin{equation}
\label{sistri1}
\left\{ \matrix{ R_\phi ( \z )  & := & - \eta^4 \partial_\phi
K_\eta (I_0 + J, \phi_0 + \widetilde \o t + \psi, \hat z) -\eta
\partial_\phi \hat G_\eta
\cr 
R_I (
\z ) & : = & \eta^2 B^T Z + \eta^4 \partial_I K_\eta (I_0 + J, \phi_0 +
\widetilde \o t + \psi, \hat z )+\eta \partial_I \hat G_{\eta}
\cr 
(R_{\overline z})_j ( \z ) & = & \ii \eta^2
(BJ)_j z_j + \ii \eta \partial_{ \bz_j} \hat G_\eta + \ii \eta^4
\partial_{ \bz_j}
K_\eta (I_0 + J, \phi_0 + \widetilde \o t + \psi, \hat z ).} \right.
\end{equation}
Remark that, since one expects $J,\psi$ and $\hat z$ to be small (they
will turn out to be of order $\eta^{2}$), and by proposition
\ref{Birknf} one has $\partial_{\overline z} \hat G_\eta(\zeta)= O(\|
\hat z\|_{a,s}^2)$, $\partial_{ \phi} \hat G_\eta(\zeta)= O(\| \hat
z\|_{a,s}^3)$, $\partial_{I} \hat G_{\eta}(\zeta)= O(\| \hat
z\|_{a,s}^3)$ the r.h.s. of (\ref{sistri},\ref{31},\ref{32}) is
actually small with respect to the period $ T \leq 2 \eta^{-2} $,
see lemma \ref{contraction}.

Define the Hilbert space $H^1_P\left((0,T);\cP_{a,s}\right)$ of the
$T$-periodic $\cP_{a,s}$ valued periodic function of class
$H^1$. In order to simplify notations we will denote this space
by $H^1_{P,s}$.
Denote 
\begin{eqnarray}
\label{norme}
|J |_{L^2,T}^2 &:=&
\frac{1}{T} \int_0^T |J|^2 \ dt\ ,\quad |\psi |_{L^2,T}^2 :=
\frac{1}{T} \int_0^T |\psi|^2 \ dt,
\\
||w||_{L^2,T,a,s}^2  &:=&   \frac{1}{T} \int_0^T ||w(t)||_{a,s}^2 \ dt
\ ,
\\
\|\zeta\|_{L^2,T,a,s}&:=& |J |_{L^2,T}+ |\psi
|_{L^2,T}+||w||_{L^2,T,a,s} \ .
\end{eqnarray}
We will endow $ H^1_P \left((0,T);\cP_{a,s}\right)\equiv H^1_{P,s}$ with
the norm
\begin{equation}
\label{norma}
\|\zeta\|_{T,a,s}:=\|\zeta\|_{L^2,T,a,s}+T\|\dot \zeta\|_{L^2,T,a,s}.
\end{equation}
\begin{remark}
\label{alg}
With this choice one has 
$$
\|\zeta(t)\|_{\cP_{a,s}}\leq C\|\zeta\|_{T,a,s}\ ,\quad \forall
t\in\real 
$$ 
with a constant independent of $T$. Therefore, with this choice of the
norm, the space $H^1_{P,s}$ is a `Banach algebra', and the $T,a,s$ norm
of the product of any component of a vector $\zeta$ with any component
of a vector $\zeta'$ is bounded by
$C\|\zeta\|_{T,a,s}\|\zeta'\|_{T,a,s} $ {\it with a constant $C$
independent of $T$}.
\end{remark}

We will consider the system (\ref{sistri},\ref{31},\ref{32}) as a
functional equation in $H^1_{P,s}$.

\begin{remark}\label{alg1}
As a consequence of (\ref{sreg}) and 
remark \ref{alg}, the map $\zeta\mapsto R(\zeta):= $ 
$(R_\phi(\zeta),$ $R_I(\zeta), $ $R_{\overline z}(\zeta))$ 
is a $C^\infty$ map from
$H^1_{P,s}$ to $H^1_{P,s+d}$. 
\end{remark}  

We are going to use the method of Lyapunov--Schmidt decomposition in
order to solve (\ref{sistri},\ref{31},\ref{32}). To this end remark
that the kernel of the linear operator $\cL$ at l.h.s. of
(\ref{sistri},\ref{31},\ref{32}) is given by $(\phi,0,0)$ with
constant $\phi\in\torus^n$. The range of $\cL$ is the space
of the functions $\zeta=(\psi,J,\hat z)$ with $\psi(t)$ having zero
mean value. So, there is a natural decomposition of $H^1_{P,s}$ in
Range+Kernel. Explicitly we write
$$
\zeta=(\phi + \psi, J , \hat z)=
(\psi, J, \hat z) + ( \phi,0,0) \equiv 
\zeta_R + \phi
$$
with $ \psi $ having zero mean value and $ \phi $ being constant.
Then we fix $\phi$, take the projection of the system
(\ref{sistri},\ref{31},\ref{32}) on the range and solve it. The
solution is a function $\zeta_R(\eta,\phi)$. Finally we insert this
function in the variational principle in order to find 
critical points of $S$.

\subsection{The Range equation}

The range equation has the form
\begin{equation}\label{sistri11}
\left\{ \matrix{
\dot J  & = & R_\phi ( \z ) - \langle R_\phi ( \z ) \rangle \label{ult1} \cr
\dot \psi - \eta^2 A \,J & = & R_I ( \z )
\cr
\dot z_j - \ii \widetilde \O_j z_j  & = & (R_{\overline z})_j ( \z ),} \right.
\end{equation}
where $\langle R_\phi ( \z ) \rangle := (1 \slash T ) \int_0^T R_\phi ( \z )
dt $. We look for its solution in the range, namely in the space 
$$
\overline H^1_{P,s} \subset H^1_{P,s}
$$ 
of the functions $\zeta_R\equiv ( 
\psi,J, \hat z ) $ with $\psi$ having zero average.

First of all we analyze the linear problem defined by the left
hand side of \equ{sistri11}. A ``small denominator problem'' appears
since inverting this linear system the denominators $ 
\widetilde \O_j T - 2 \pi l $, $j \geq n+1$, 
$ l \in \integer $ are present. So, define the linear operator
$$
\cL(\psi,J,\hat z)\equiv \cL\zeta_R:=(\dot J,\dot\psi-\eta^2AJ,\dot
w_j- \ii \widetilde \Omega_jw_j)
$$
and study the equation
\begin{equation}
\label{lineariz}
\cL\zeta_R= (\widetilde \psi,\widetilde J,\widetilde w)
\end{equation}
with $(\widetilde \psi,\widetilde J,\widetilde w)\in \overline 
H^1_{P,s+ \tau}$ given.

\lem{lemma 1}
Assume (H2). If 
$$
\widetilde \zeta_R\equiv ( \widetilde \psi,\widetilde J, \widetilde w) \in
\overline H^1_{P,s+\tau}\quad {\rm i.e.\ in}\ H^1_{P,s + \tau} 
\ {\rm with } \ \
\int_0^T {\widetilde \psi }(t) \ dt = 0, 
\label{spaceY}
$$
then the equation (\ref{lineariz}) has a unique solution
$$ 
\zeta_R \equiv (\psi,J, w ) \in \overline H^1_{P,s}.  
$$
Moreover, for $ T \in (\eta^{-2}, 2 \eta^{-2} )$ and
a constant $ C:= C(\delta )$   
$$
\|\zeta_R\|_{T,a,s}\leq \frac{C}{\eta^2}\|\widetilde \zeta_R\|_{T,a,s+\tau}.
$$ 
\elem

\noindent {\bf Proof.}
Since $A$ is symmetric and invertible it has an orthonormal basis
of eigenvectors $ e_1, \ldots , e_n $ with eigenvalues 
$\lambda_1, \ldots, \lambda_n $. In these coordinates 
$ J(t) = \sum_{k=1}^n J_k(t) e_k$,  $ \psi(t) = \sum_{k=1}^n \psi_k (t) e_k$ 
and the solution $ \zeta_R $ of (\ref{lineariz}) with $ \psi_0 = 0 $ 
has Fourier coefficients
$$
J_{kl} = \frac{T {\widetilde \psi}_{kl} }{ \ii 2 \pi l } \quad {\rm for} \quad
l \neq 0, \qquad  J_{k0} = - \frac{ {\widetilde J}_{k0}}{\eta^2 \lambda_k}, 
$$
$$
\psi_{kl} = T \frac{ {\widetilde J}_{kl} + \eta^2 J_{kl} 
\lambda_k}{\ii 2 \pi l } \quad {\rm for} \quad
l \neq 0, 
$$
and, for $ j \geq n + 1 $ 
$$
w_{jl} := \frac{T {\widetilde w}_{jl}}{\ii ( 2 \pi l - \widetilde{\O}_j T)}.
$$
We find
\beq{es1}
|J|_{L^2,T}^2 = \sum_{kl} |J|_{kl}^2 = \sum_k 
\Big( \frac{ {\widetilde J}_{k0}}{\eta^2 \lambda_k} \Big)^2 +
\sum_{k,l \neq 0 } \Big( 
\frac{T {\widetilde \psi}_{kl} }{\ii 2 \pi l } \Big)^2 \leq
\frac{C}{\eta^4} | \widetilde J |_{L^2,T}^2 + C T^2 | \widetilde 
\psi |_{L^2,T}^2.
\eeq
A similar estimate for $ | \psi |_{L^2,T} $ holds. Moreover
\beq{es2}
| { \dot \psi }|_{L^2,T} \leq \eta^2 | J|_{L^2,T} + 
|\widetilde \psi|_{L^2,T} \leq C ( | \widetilde J|_{L^2,T} 
+ |\widetilde \psi |_{L^2,T}),
\eeq
using \equ{es1}.
Finally, 
the solution $ w = ( w_j  )_{j\geq n+1} $ 
of (\ref{lineariz}) is 
$$
w_j ( t ) = \sum_{ l \in \integer } 
\frac{T {\widetilde w}_{jl}}{\ii ( 2 \pi l -
\widetilde{\O}_j T)} e^{\ii (2 \pi \slash T) l t} 
$$
where $ {\widetilde w}_j ( t ) = $  
$ \sum_{ l \in \integer } {\widetilde w}_{jl} 
e^{\ii (2 \pi \slash T) l t}. $
From (H2) we get
\beq{stw}
|| w ||_{L^2,T,a,s} \leq C \frac{T}{\d} || {\widetilde w}||_{L^2,T,
a, s + \tau},
\qquad 
|| {\dot w} ||_{L^2,T,a,s} 
\leq C \frac{T}{\d} || \dot {\widetilde w}||_{L^2,T,a,s +\tau}.
\eeq
By \equ{es1}, \equ{es2} and \equ{stw} the last estimate of the lemma 
follows. 
\hfill$\Box$

\bigskip

Thus $\cL^{-1}$ defines a linear bounded operator 
$L:\overline H^1_{P,s+\tau}\to \overline H^1_{P,s}$. 

\giu

In order to find a solution $ \z_R = (\psi,J, \hat z ) $ 
of the range equation it is sufficient to find a fixed point of 
\beq{pfisso} 
\z_R = \Phi ( \z_R ) := 
L \big( N( \z_R ; \phi ) \big) 
\eeq 
in the space $ \overline
H^1_{P,s}$, where $ N := N ( \z_R ; \phi)
$ denotes the right hand side of \equ{sistri11}. 

\begin{lemma}
\label{contraction} Assume $ d > \tau $. Then
there exists a constant $ C $ sufficiently large 
such that $ \forall \eta \ll 1 $ the map
$\Phi$ is a contraction of a ball of radius $ C \eta^2 $.
\end{lemma}

\noindent{\bf Proof}. Consider a $\zeta_R\in\overline H^1_{P,s}$ with
$\|\zeta_R\|_{T,a,s}\leq \rho$ with some positive (small)
$\rho$. Since $H_{P,s}^1$ is an algebra with constants independent of
$T$ (c.f. remark \ref{alg}), one has, by (\ref{sistri1}),
$$
\|N(\zeta_R)\|_{T,a,s+d}\leq C(\eta^4+\eta\rho^2)
$$ 
with a suitable $ C $. Therefore, by lemma \ref{lemma 1} one has
$$
\|\Phi(\zeta_R)\|_{T,a,s}\leq\|\Phi(\zeta_R)\|_{T,a,s+d-\tau} \leq
C\left(\eta^2+\frac{\rho^2}{\eta}\right) 
$$ which is smaller than $\rho$ provided $C(\eta^2+\rho^2/\eta)<\rho$,
which is implied e.g. by $\rho = 2 C \eta^{2} $ and $\eta$ 
small enough.

Similarly one estimates the Lipschitz constant of $\Phi$ by the norm
of its differential. Such a differential is bounded in a ball of
radius $\rho$ by $C(\eta^2+\rho/\eta)$, from which the thesis
follows. \hfill $\Box$

\begin{corollary}
\label{range}
There exists a unique smooth function $\torus^n\ni\phi\mapsto
\zeta_R(\phi,\eta)\in\overline H^1_{P,s}$ solving (\ref{sistri11}) and
fulfilling
$$
\|\zeta_R(\phi,\eta)\|_{T,a,s}\leq C\eta^{2} \ .
$$
\end{corollary}

\subsection{The kernel equation}

The geometric interpretation of the construction of the previous
subsection is that we have found a submanifold $\cT^n\equiv
(\phi,\zeta_R(\phi,\eta))\subset H^1_{P,s}$, diffeomorphic to an $n$
dimensional torus, on which the partial derivative of the action
functional $S$ with respect to the variables $ \zeta_R $
vanishes. Consider now the restriction $ S_{n} $ of $S$ to $\cT^n$. 
At a critical point of $S_n$ all the partial derivatives
of the complete functional $S$ vanish, and therefore we can conclude
that such point is critical also for the non restricted functional. 

By standard Lusternik-Schnirelmann theory there exist at least $n$
geometrically distinct $T$-periodic solutions, i.e. solutions not
obtained one from each other simply by time-translations.  Indeed,
restrict $ S_{n} $ to the plane $ E := [\widetilde \o ]^{\bot} $
orthogonal to the periodic flow $ \widetilde \o = (1/T) 2 \pi k $ with
$ k \in \integer^n $.  The set $ \integer^n \cap E $ is a lattice of $
E $, and hence $ S_{n} $ defines a functional $S_{n|\Gamma}$ on the
quotient space $\Gamma := E \slash (\integer^n \cap E) \sim
\torus^{n-1}$.

Due to the invariance of $S_n$ with respect to the time shift, a
critical point of $ S_{n|\Gamma}$ is also a critical point of $ S_n :
\torus^n \to \real $.  By the Lusternik-Schnirelman category theory
since cat$\Gamma$ = cat$\torus^{n-1}$ = $n$, we can define the $n$
min-max critical values $c_1 \leq c_2 \leq \ldots \leq c_n $ for $
{S_n}_{| \Gamma } $.  If the critical levels $ c_i $ are all distinct,
the corresponding $T$-periodic solutions are geometrically distinct,
since their actions $ c_i $ are all different. On the other hand, if
some min-max critical level $ c_i $ coincide, then, by the
Lusternik-Schnirelmann 
theory, $ {S_n}_{| \Gamma } $ possesses
infinitely many critical points. However {\it not} all the
corresponding $T$-periodic solutions are necessarily geometrically
distinct, since two different critical points could belong to the same
orbit. Nevertheless, since a periodic solution can cross $ \Gamma $ at
most a finite number of times, the existence of infinitely many
geometrically distinct orbits follows. For further details see [BBV03].

This concludes the proof of the theorem \ref{main?}.

\section{Applications}

\subsection{The nonlinear beam equation}

Consider the beam equation
\beq{maineq}
u_{tt} + u_{xxxx} + mu = f ( u ) 
\eeq
subject to hinged boundary conditions
\beq{BC}
u(0,t)= u(\pi,t)= u_{xx}(0,t)= u_{xx}(\pi,t)=0,
\eeq
where the nonlinearity $ f (u) $ is a real analytic, odd function 
of the form 
$$
f(u) = a u^3 + \sum_{k \geq 5} f_k u^k, \qquad a \neq 0.
$$
The beam equation \equ{maineq} is a Hamiltonian PDE
with associated Hamiltonian 
$$
H = \int_0^\pi \frac{u_t^2}{2} + \frac{ u_{xx}^2}{2}  
+ \frac{m u^2}{2} + g(u) \ dx
$$
where $ g(u) := \int_0^u f(s) \ ds $ is a primitive of $ f $. 

Write the system in first order form
\begin{equation}
\label{1beam}
\left\{ \matrix{\dot u&=&v
\cr
\dot v&=& - u_{xxxx} - mu +f ( u ) 
}\right.
\end{equation}
The standard phase space\footnote{An equivalent definition makes use
of the so called compatibility conditions required for the smoothness
of solutions of second order equations with Dirichlet boundary
conditions, see e.g. [Bre93], theorem X.8.} for (\ref{1beam}) is
$\cF_s:=H^{s}_C\times H^{s-2}_C\ni(u,v)$, where $H^s_C$ is the space
of the functions which extend to skew symmetric $H^s$ 
periodic functions over $[-\pi,\pi]$. Note that 
$ H^s_C = \{ u(x) = \sum_{j \geq 1} u_j \sin (jx) \ | \
\sum_{j \geq 1} |u_j |^2 j^{2s} < + \infty \} $.
It is then immediate to realize that, due to the
regularity and the skewsymmetry of the vector field of the nonlinear
part, $f$ defines a smoothing operator, namely a smooth map from $\cF_s$
to $\cF_{s+2}$, provided $s\geq 1$.

Here we are also interested in spaces of analytic functions, namely
functions whose Fourier coefficients belong to $\cH^{a,s}$ with some
positive $a$. It is easy to see that the smoothing
property of the nonlinearity holds also for these spaces.

Introduce coordinates 
$q = (q_1, q_2, \ldots, )$, $p = (p_1, p_2, \ldots, )$
through the relations 
$$
u(x) = \sum_{j \geq 1} \frac{q_j}{\sqrt{\o_j}} \phi_j (x), \qquad
v(x) = \sum_{j \geq 1} p_j \sqrt{\o_j} \phi_j (x),
$$
where $\phi_j (x) = \sqrt{2 / \pi} \sin (j x) $ and
\begin{equation}
\label{j4}
\o_j^2 = j^4 + m. 
\end{equation} 
Remark also that 
$$
\omega_j\sim j^2\ .
$$
Passing to complex coordinates 
$$
z_j   := \frac{q_j + \ii p_j }{ \sqrt{2} }, \qquad 
\bz_j := \frac{q_j - \ii p_j }{ \sqrt{2} },
$$
the Hamiltonian takes the form (\ref{hami}) and the nonlinearity
fulfills (S) with $s\geq 1$ a suitable $a$, depending on the
anayticity strip of $f$, and $d=2$ (for more details see 
[P\"os96b]-[GY03]).

In order to verify the nonresonance property we use $m$ as a parameter
belonging to the set $[0,L]$ with an arbitrary $L$. 

\begin{lemma}
\label{reso}
There exists a finite set $\Delta\subset [0,L]$ such that, if
$m\in[0,L]\string\ \Delta$ then condition
(NR) holds.
\end{lemma}

\noindent {\bf Proof.} First remark that, due to the growth property
of the frequencies there is at most a finite number of vectors
$l\in\integer^2$ at which $\omega\cdot k+\Omega\cdot l$ is small. It
follows that, having fixed an arbitrary constant $C$, there is at most
a finite set of $k$'s and $l$'s over which $|\omega\cdot k+\Omega\cdot
l|<C$. Denote by $\cS$ such a set.

For $(k,l)\in\cS$ consider 
$$
f_{kl}(m)=\omega(m)\cdot k+\Omega(m)\cdot l 
$$
since $f_{kl}$  is an analytic function it has only isolated zeros. So at
most finitely many of them fall in $[0,L]$. The set $\Delta$ is the
union over $k,l\in\cS$ of such points. Fix $m\in [0,L]\backslash
\Delta$. \hfill$\Box$

\bigskip

Then one can put the system in seminormal form. The explicit
computation was essentially done in [KP96] (see also [P\"os96b,GY03]), obtaining
that the matrixes $A$ and $B$ are given by
\begin{equation}
A = \frac{6}{\pi}
\left( \matrix{ \frac{3}{\o_1^2} && \frac{4}{\o_1 \o_2} && \ldots  && 
\frac{4}{\o_1 \o_n} \cr
\frac{4}{\o_2 \o_1} && \frac{3}{\o_2^2} && \ldots  && 
\frac{4}{\o_2 \o_n} \cr
\ldots && \ldots && \ldots  && \ldots \cr
\frac{4}{\o_1 \o_n} && \frac{4}{\o_n \o_2} && \ldots  && 
\frac{3}{\o_n^2} 
} \right),
\ 
B = \frac{6}{\pi}
\left( \matrix{ \frac{4}{\o_{n+1} \o_1} && \ldots  && 
\frac{4}{\o_{n+1} \o_n} \cr
\frac{4}{\o_{n+2} \o_1} && \ldots  && 
\frac{4}{\o_{n+2} \o_n} \cr
\vdots && \vdots  && \vdots \cr
} \right).
\end{equation}
Remark that, defining the matrixes $S_1:={\rm diag} (\omega_1, ... , \omega_n)$
and $S_2:={\rm diag} (\omega_{n+1}, \omega_{n+2}, ...)$ one can write $
A=\frac{6}{\pi}S_1^{-1}\widetilde AS_1^{-1}\ ,\quad B=\frac{6}{\pi}
S_2^{-1}\widetilde BS_1^{-1} $ with
\begin{equation}
\label{tildeA}
\widetilde A = 
\left( \matrix{ {3} && {4}&& \ldots  && 
{4}\cr
{4} && {3} && \ldots  && 
{4} \cr
\ldots && \ldots && \ldots  && \ldots \cr
{4} && {4} && \ldots  && 
{3} 
} \right),
\ 
\widetilde B =
\left( \matrix{ {4} && \ldots  && {4} \cr
{4} && \ldots  && 
{4} \cr
\vdots && \vdots  && \vdots \cr
} \right).
\end{equation}
With these expressions at hand it is immediate to verify that
det$A\not=0$. For what pertains $\hat \Omega_j$ by exactly the same
argument in the proof of lemma \ref{reso} one has that they are
different from zero except for at most finitely many values of $m\in[0,L]$.

{\it Thus, provided $m$ does not belong to a finite subset of $[0,L]$,
theorem \ref{main?} and its corollary \ref{main} apply. }

\subsection{A nonlinear Schr\"odinger equation}

Consider the space $H^s_C$ as in the previous section.
Following P\"oschel [P\"os02] we define a smoothing operator as follows.
Fix a sequence $\{\rho_j\}_{j\geq1}$ with the property 
\begin{equation}
\label{schro0}
\forall j\geq 1\ \quad
\rho_j\not=0\ ,\ {\rm and}\quad
|\rho_j|\leq Cj^{-d/2}\ ,\quad d>1.
\end{equation}
Consider the even, $2\pi$ periodic function 
$ \rho(x) := \sum_j \rho_j \cos(jx) $ and define
\begin{equation}
\label{schro1}
\Gamma: H^s_C\to H^{s+d/2}_C \ , \quad \Gamma u := \rho * u
\end{equation}
where the star denotes convolution (it is defined 
first extending the function $ u $ to an odd $2\pi$ periodic function).

\begin{remark}
\label{smoo}
It is easy to see that, expanding $u$ in Fourier series
$$
u(x)=\sum_{j\geq 1}{z_j}\sin(jx), 
$$
the $j$-th Fourier coefficient of $\Gamma u$ is proportional to
$\rho_jz_j$.
\end{remark}

Consider the Hamiltonian system with Hamiltonian function 
\begin{equation}
\label{schro2}
H(u,\overline u)=\int_0^\pi |u_x|^2+
P\left(\left|\Gamma u\right|^2\right)
\end{equation}
with $P$ an analytic function having a zero of order two at
the origin. So the equations of motion are
\begin{equation}
\label{shro3}
- \ii \dot u=u_{xx}+
 \Gamma\left(P'\left(\left|\Gamma u\right|^2\right)\Gamma u\right).
\end{equation}
Inserting the Fourier expansion of $u$ the Hamiltonian takes the form
 (\ref{hami}) with $\omega_j=j^2$ and the vector field fulfills (S)
 with $d$ given by (\ref{schro0}).

Then one has that (NR) is here violated. So in principle one cannot
expect proposition \ref{Birknf} to hold for this system. However it is
clear that one can use a similar procedure to obtain a weaker result
in which the function $\overline G$ is in resonant normal form, i.e. it
contains only monomyals Poisson commuting with $H_0$. Now, an explicit
computation, identical to that done in [KP96], shows that the so obtained
normal form {\it depends on the actions only}. So one can proceed in
the construction and try to look for resonant tori and the
corresponding periodic orbits. 

The explicit computation shows that also in this case the matrices $A$
and $B$ have the structure $ A=\alpha S_1\widetilde AS_1\,$, 
$B = 2\alpha S_2\widetilde BS_1 $ with the matrixes $\widetilde A$ and
$\widetilde B$ still given by (\ref{tildeA}), matrixes $S_1:={\rm
diag}(\rho_1, ... , \rho_n)$ and $ S_2 := {\rm diag} (\rho_{n+1},
\rho_{n+2}, .... )$ and a suitable constant $ \alpha $. So the
determinant of $A$ is still different from zero. The frequencies $\hat
\Omega_j$ have now the structure
$$
\hat\Omega_j(\rho)=j^2-\rho_j a(\rho)\ ,\quad\forall j\geq n+1
$$
where $a$ is a function of $\rho_1,...,\rho_n$.  So, except for
exceptional choices of $\{\rho_j\}_{j\geq n+1}$ the nondegeneracy
conditions are fulfilled and the conclusions of 
theorem \ref{main?} and its corollary apply.

\section{Appendix: Proof of proposition \ref{Birknf}}

The idea is to proceed as in the proof of the standard
Birkhoff normal form theorem, i.e. by successive elimination of the
nonresonant monomials. As a variant with respect to the 
standard procedure 
one does not eliminate terms which are at least cubic in the 
variables $ \hat z$. Remark that the estimates involved in the proofs
are much more complicated than in the finite dimensional case. 

To start with expand  $P$ in Taylor series up to order five: $
P=P_3+P_4+P_5+$higher order terms.  Then we begin by looking for the
transformation simplifying $P_3$. So write
\begin{eqnarray*}
P_3&=&P_3^1+O(\|\hat z\|^3) \ ,
\end{eqnarray*}
where $P_3^1$ is composed by the first three terms of the Taylor
expansion of $P_3$ in the variables $\hat z$ only (so it contains only
terms of degree 0,1 and 2 in such variables).  We use the Lie
transform to eliminate from $P_3^1$ all the nonresonant terms, i.e. we
make a canonical transformation which is the time 1 flow $\Phi^1$ of
an auxiliary Hamiltonian system with a Hamiltonian function $\chi$ of
degree 3. By considering the Taylor expansion of $\Phi^1$ at zero one
has
\begin{equation}
\label{norpro2}
H\circ \Phi^1= H_0+P_3^1+\left\{\chi,H_0 \right\}+O(\| z\|^4)
+O(\|\hat z\|^3) \ .
\end{equation}
One wants to determine $\chi$ so that 
$P_3^1+\left\{\chi,H_0 \right\}$ depends on the actions $|z_j|^2$
only. To this end we proceed as usual in the theory of Birkhoff normal
form. 

Denote by $x=(x_1,...,x_n)\equiv (z_1,..., z_n)$ the first $n$
variables and take $\chi$ to be a homogeneous polynomial of degree
three. Write
\begin{equation}
\label{norfor4}
\chi=\sum_{|j_1|+|j_2|+|j_3|+|j_4|=3} \chi_{j_1j_2j_3j_4} x^{j_1}\overline
x^{j_2}\hat z^{j_3}\overline {\hat z}^{j_4}
\end{equation}
with multindexes $j_1,j_2,j_3,j_4$.  For a multi index $j_l\equiv
(j_{l,1},...,j_{l_n})$ we used the notation $| j_l | := |j_{l,1} | +
\ldots |j_{l,n} | $ and $x^{j_l} := x_1^{j_{l,1}} \ldots x_n^{j_{l,n}}
$, and similarly for a multi index with infinitely many components.
So, one has
$$ \left\{\chi,H_0\right\}=
\sum_{|j_1|+|j_2|+|j_3|+|j_4|=3}\ii\left(\omega\cdot(j_1-j_2)+\Omega\cdot
(j_3-j_4)\right) \chi_{j_1j_2j_3j_4} x^{j_1}\overline x^{j_2}\hat
z^{j_3}\overline {\hat z}^{j_4}.
$$
Write now 
\begin{equation}
\label{norfor5}
P_3^1=\sum_{|j_1|+|j_2|+|j_3|+|j_4|=3} P_{j_1j_2j_3j_4} x^{j_1}\overline
x^{j_2}\hat z^{j_3}\overline {\hat z}^{j_4}
\end{equation}
and remark that the indexes are here subjected to the further limitation
$|j_3|+|j_4|\leq 2$. So one is led to the choice
\begin{equation}
\label{chi}
\chi_{j_1j_2j_3j_4} :=\frac{-P_{j_1j_2j_3j_4}}{ \ii
  \left(\omega\cdot(j_1-j_2)+\Omega\cdot 
(j_3-j_4)\right) } \ ,\quad j_1-j_2+j_3-j_4\not=0
\end{equation}
and zero otherwise.  Remark that, due to the assumption (NR) the
denominators appearing in the above expression are all different from
zero. Moreover, since due to the growth of the frequencies they are
actually bounded away from zero.  Then in order to conclude the proof
(at least for what concerns the elimination of the third order part)
one has to ensure that the function $\chi$ is well defined and that it
has a smooth Hamiltonian vector field. The terms of $\chi$ of
different degree in $\hat z$ have to be treated in a different way, so
we will denote by $\chi_0,\chi_1,\chi_2$ the homogeneous parts of degree
0,1 and 2 respectively with respect to the variables $\hat z$.

We need a few lemmas. 

\begin{lemma}
\label{app1}
Let $\real^n\ni x\mapsto f(x)\in\ell^2$ be a homogeneous bounded
polynomial of degree $r$. Write
$$
f(x)=\sum_{j\in\natural^n,|j|=r}\sum_{k\geq1}f_{jk}x^je_k
$$
where $e_k$ is the standard basis of $\ell^2$. Let
 $\left\{\rho_{j,k}\right\}_{j\in\natural^n}^{k\geq 1}$
be a sequence with the property $|\rho_{jk}|\geq C$ and define a
 function $g$ by
\begin{equation}
\label{app25}
g(x)=\sum_{j\in\natural^n,|j|=
r}\sum_{k\geq1}\frac{f_{jk}}{\rho_{jk}}x^je_k.
\end{equation}
Then there exists $C$ such that $ \|g(x)\|\leq C\|x\|^r.$
\end{lemma}
\noindent{\bf Proof.} Write $ g(x)=\sum_{j}g_j x^j$ 
and remark that the cardinality of the set over which the sum is
carried out is finite. We estimate each of the vectors $g_j$'s. One
has
$$
\|g_j\|^2 =\sum_k \left(\frac{f_{kj}}{\rho_{kj}}\right)^2\leq
\frac{1}{C^2}\sum_{k}f_{kj}^2=  \frac{1}{C^2}\|f_j\|^2.
$$
Now the norms of the vectors $f_j$ are bounded and therefore the
thesis follows.\hfill $\Box$

\begin{remark}
\label{app32}
By the same proof, the same result holds if the space
$\ell^2$ is substituted by the spaces $\cH^{a,s}$.
\end{remark}

\begin{lemma}
\label{app2}
Let $\real^n\times \ell^2\ni (x,z)\mapsto f(x,z)\in\real$ be a
homogeneous bounded polynomial of degree $r$ in $x$, linear and bounded
in $z$. Write
$$
f(x,z)=\sum_{{k\geq1 \atop j\in\natural^n,|j|=r}}f_{jk}x^jz_k.
$$
Let
$\left\{\rho_{j,k}\right\}_{j\in\natural^n}^{k\geq 1}$ be as above
and define a function $g$ by
\begin{equation}
\label{app22}
g(x,z)=\sum_{{k\geq1 \atop j\in\natural^n,|j|=r}}
\frac{f_{jk}}{\rho_{jk}}x^jz_k.
\end{equation}
Then there exists $C$ such that 
$ |g(x,z)|\leq C\|x\|^r\|z\|.$
\end{lemma}
\noindent{\bf Proof.} Just write 
$ g(x,z)=\sum_{j}g_j(z) x^j\ .$ 
Fix $j$ and study  the linear functional $g_j(z)$, one has
$$
|g_j(z)|=\left|\sum_{k\geq 1}f_{jk}\frac{z_k}{\rho_{jk}} \right| \leq
\|f_j\| \|z/\rho\|
$$ where $f_j$ is defines in analogy to $g_j$, its norm is the norm as
a linear functional, and $z/\rho$ is the vector of $\ell^2$ with $k$-th
component equal to $z_k/\rho_{jk}$. From this inequality, summing over
$j$, the thesis follows. \hfill$\Box$

\bigskip

In order to estimate the vector field of $\chi_2$ we will need the
following lemma:
\begin{lemma}
\label{poschel}
[Lemma A.1 of [P\"os96a]]. If $A=(A_{kl})$ is a bounded linear operator on
$\ell^2$, then also $B=(B_{kl})$ with
\begin{equation}
\label{posch2}
B_{kl}:=\frac{|A_{kl}|}{1+|k-l|}\ ,
\end{equation} 
is a bounded linear operator on $\ell^2$.
\end{lemma}
For the proof we refer to [P\"os96a].

\begin{lemma}
\label{app3}
Let $\real^n\times \ell^2\ni (x,z)\mapsto f(x,z)\in\ell^2$ be a
homogeneous bounded polynomial of degree $r$ in $x$ linear and bounded
in $z$. Write
$$
f(x,z)=\sum_{{k,l\geq1 \atop j\in\natural^n,|j|=r}}f_{jkl}x^jz_ke_l.
$$
Let
 $\left\{\rho_{j,k,l}\right\}_{j\in\natural^n}^{k,l\geq 1}$ be a sequence
fulfilling
\begin{equation}
\label{sotto}
|\rho_{jkl}|\geq C_1(1+|k-l|)
\end{equation}
and define a
 function $g$ by
\begin{equation}
\label{app24}
g(x,z)=\sum_{{k,l\geq1 \atop j\in\natural^n,|j|=r}}
\frac{f_{jkl}}{\rho_{jkl}}x^jz_k.
\end{equation}
Then there exists $C$ such that $ \|g(x,z)\|\leq C\|x\|^r\|z\|.$
\end{lemma}
\noindent{\bf Proof.}  Write 
$ g(x,z)=\sum_{j}g_j(z) x^j\ .$ 
Fix $j$ and apply lemma
\ref{poschel} to such operators obtaining the result. \hfill$\Box$

\begin{remark}
\label{pos2}
An identical statement holds for functions from $\real\times
\cH_{a,s}$ to $\cH_{a,s+d}$. To obtain the proof just remark that the
boundedness of a linear operator $B=(B_{kl})$ ($g_j$ in the proof) as an operator
from $\cH^{a,s}$ to $\cH^{a,s+d}$ is equivalent to the boundedness of
$\tilde B:=(v_k B_{kl} s_l)$ as an operator from $\ell^2$ to
itself, where $v_k,s_l$ are suitable weights.
\end{remark}

\bigskip

With the above lemmas at hand it easy to estimate the vector field of
$\chi$. We treat explictly only $\chi_1$. 
\begin{lemma}
\label{chi1}
Let $\chi_1$ be the component linear in $\hat z$ and
$\overline{\hat z}$ of the function $\chi$ defined by (\ref{chi}). Then
there exists a constant $C$ such that its vector field is bounded by
$$
\|X_{\chi_1}(z,\overline z )\|_{a,s+d}\leq C \|z\|_{a,s}^2.   
$$
\end{lemma}
\noindent {\bf  Proof.} Write $\chi_1$ as follows
$$
\chi_1(x,\overline x,\hat z,\overline{\hat z})= \left\langle \chi_{01}(x,\overline
x); \hat z \right\rangle_{\ell^2}+ \left\langle \chi_{10}(x,\overline
x); \overline {\hat z} \right\rangle_{\ell^2}.
$$
Consider the first term. Separating the $x,\overline x$ and $\overline {\hat z}$
components, its vector field is given by  
$$ \left(\ii \left\langle \frac{\partial \chi_{01}}{\partial \overline x};
\hat z \right\rangle_{\ell^2} , - \ii \left\langle \frac{\partial
\chi_{01}}{\partial \overline x}; \hat z \right\rangle_{\ell^2},-\ii
\chi_{01}(x,\overline x) \right).
$$
Explicitly $\chi_{01}$ is given by 
$$
\sum_{{|j_1|+|j_2|=2\atop l\geq
    n+1}}\frac{- P_{j_1j_2e_l}}{\ii(\omega\cdot
  (j_1-j_2)+\Omega_l)}x^{j_1}\overline x^{j_2} e_l.
$$
It follows that each of the $x$ (and $\overline x$) components of the
vector field has the structure considered in lemma \ref{app2}, which
therefore gives the estimate of such part of the vector
field. Concerning the $\overline {\hat z}$ component lemma \ref{app1}
applies and gives the result. The remaining components can be treated
exactly in the same way.\hfill$\Box$

\bigskip

The estimate of the vector fields of $\chi_0$ and $\chi_2$ are
obtained in a similar way. In order to apply lemma \ref{app3} to the
estimate of the vector field of $\chi_2$ one has just to remark that
from (A) and (NR) one has the
estimate
$$
|\omega\cdot k+\Omega_j-\Omega_l|\geq C(1+|j-l|)\ . 
$$
Thus we have the following 
\begin{proposition}
\label{campo}
The vector field of the function $\chi$ defined by (\ref{chi})
fulfills the inequality
$$
\|X_{\chi}(z,\overline z)\|_{a,s+d}\leq C\|z\|_{a,s}^2\ .
$$
\end{proposition}

Then by standard existence and uniqueness theory one has that such
vector fields defines a unique smooth time 1 flow in a neighbourhood
of the origin both in $\cP_{a,s}$ and in $\cP_{a,s+d}$. It follows
that the transformation is well defined. Transforming the vector
field of $H$ one gets a vector field having the same smoothness
properties of the original one. Thus one can iterate the construction
and eliminate also the unwanted terms of degree four and five. This
concludes the proof of proposition \ref{Birknf}.

\section*{References}

\footnotesize

\noindent[ACE87]
Ambrosetti, A.; Coti-Zelati, V., Ekeland, I.:
{\sl Symmetry breaking in Hamiltonian systems},
Journal Diff. Equat. 67, 1987, p. 165-184.
\smallskip

\noindent[Bam00] Bambusi, D.: {\sl 
Lyapunov Center Theorems for some nonlinear
PDE's: a simple proof}, Ann. Sc. Norm. Sup. di Pisa, Ser. IV, 
{\bf XXIX}, 2000.
\smallskip

\noindent[Bam03] Bambusi, D.: {\sl Birkhoff normal form for some
nonlinear {PDE}s}, Commun. Math. Phys., {\bf 234}, (2003)
	 {253-285}. 
\smallskip

\noindent[BG03] Bambusi D., Grebert B.: {\sl Forme normale pour NLS en
dimension quelconque}, C.R. Acad. Sci. Paris Ser. 1, {\bf 337}, (2003)
409--414. 
\smallskip

\noindent[BP01] Bambusi D., Paleari S.: {\sl Families of periodic
orbits for resonant {PDE}'s}, {J. Nonlinear Science}, {\ bf 11},
({2001}), {69--87}.

\smallskip
\noindent[BBV03]
Berti, M., Biasco, L., Valdinoci, E.: {\sl Periodic orbits
close to elliptic tori and applications to the three body problem},
preprint 2003, available at
{{\tt http://www.math.utexas.edu/mp\_arc.
}}
\smallskip

\noindent[BB03] Berti M., Bolle P.:
{\sl Periodic solutions of nonlinear wave equations with general 
nonlinearities}, Commun. Math. Phys. to appear. 

\smallskip

\noindent[BL34]
Birkhoff, G.D., Lewis, D.C.:
{\sl On the periodic motions near a given periodic
motion of a dynamical system},
Ann. Mat. Pura Appl., {\bf IV}. Ser. 12, 1934, 117-133.
\smallskip

\noindent [Bou95] Bourgain J.:
 {\sl Construction of periodic solutions 
of nonlinear wave equations in higher dimension}, {Geometric and
 Functional Analysis}, {\bf 5}, ({1995}) 629-639.

\smallskip


\noindent[Bre93] Brezis, H.: {\it Analyse fonctionelle}, Masson, Paris
1993. 

\smallskip
\noindent [CW93] Craig W., Wayne C.E.: {\sl Newton's method and
periodic solutions of nonlinear wave equations}, {Comm. Pure
  Appl. Math.}, {\bf 46} {1993},
{1409-1498}.

\smallskip
\noindent[GY03]
Geng J., You, J.:
{\sl KAM tori of Hamiltonian perturbations of 1D linear beam equations},
J. Math. Anal. Appl. {\bf 277}, 2003, 104-121.
\smallskip

\noindent [Kuk93] Kuksin S.B.: {Nearly integrable infinite-dimensional
{H}amiltonian Systems}, Springer-Verlag, {Berlin}, 1993.

\smallskip
\noindent[KP96] Kuksin S.B., P\"oschel J.: {\sl Invariant {C}antor
manifolds of quasi-periodic oscillations for a nonlinear
{S}chr\"odinger equation}, {Ann. of Math.}, {\bf 143}, {1996},
{149-179}.
\smallskip

\noindent[Lew34] Lewis, D.C.: {\sl Sulle oscillazioni periodiche di un
sistema dinamico}, Atti Acc. Naz. Lincei, Rend. Cl. Sci. Fis. Mat.
Nat., {\bf 19}, 1934, 234-237.
\smallskip

\noindent[Mos77] Moser, J.:{ \sl Proof
of a generalized form of a fixed point Theorem due to G. D. Birkhoff},
Geometry and topology (Proc.
III Latin Amer. School of Math.,
Inst. Mat. Pura Aplicada CNP, Rio de Janeiro, 1976),
pp. 464--494. Lecture Notes in Math.,
Vol. 597, Springer, Berlin, 1977.
\smallskip

\noindent[P\"os96a]  P\"oschel, J.: {\sl A KAM-theorem for some nonlinear
partial differential equations}, Ann. Sc. Norm. Sup. di Pisa, Ser. IV, 
{\bf XXIII}, (1996), 119-148. 
\smallskip

\noindent[P\"os96b]  P\"oschel, J.: {\sl Quasi-periodic
solutions for a nonlinear wave equation}, Comment. Math. Helv. {\bf 71}
(1996), 269--296. 

\smallskip

\noindent[P\"os02]  P\"oschel, J.:{ \sl  On the construction of almost
periodic solutions for a nonlinear Schr\"odinger equation}, Ergodic
Theory Dynam. Systems {\bf 22} (2002), 1537--1549.

\smallskip
\noindent[Way90] Wayne C.E.:
{\sl Periodic and quasi-periodic solutions of nonlinear
wave equations via {KAM} theory}, {Commun. Math. Physics}, {\bf 127}
({1990}) {479-528}.

\end{document}